\documentclass[12pt]{amsart}

\usepackage{amsfonts,amssymb,amsthm,latexsym}

\newcommand{\TDS}{\renewcommand{\baselinestretch}{1.7}\tiny\normalsize}

\newtheorem{thm}{Theorem}[section]
\newtheorem{prop}[thm]{Proposition}

\newtheorem{cor}[thm]{Corollary}
\newtheorem{con}[thm]{Conjecture}

\newtheorem{claim}[thm]{Claim}

\newcommand{\xni}{$ \, ({{x}_n}_i) \,$ }
\newcommand{\xn}{$ \,(x_n)_{n=1}^\infty \,$ }

\newcommand{\lb}{\label}
\newcommand{\bi}{\bibitem}
\newcommand{\ci}{\cite}
\newcommand{\no}{\nonumber}

\newcommand{\be}{\begin{equation}}
\newcommand{\ee}{\end{equation}}

\newcommand{\bea}{\begin{eqnarray}}
\newcommand{\eea}{\end{eqnarray}}
\newcommand{\bean}{\begin{eqnarray*}}
\newcommand{\eean}{\end{eqnarray*}}
\newcommand{\ba}{\begin{array}}
\newcommand{\ea}{\end{array}}
\newcommand{\mb}{\mbox}

\newcommand{\bd}{\begin{displaymath}}
\newcommand{\ed}{\end{displaymath}}

\begin{document}

\title{Genus n  Banach Spaces}
\author{P. G. Casazza  and M. C. Lammers}
\address{Department of Mathematics \\
University of Missouri-Columbia \\
Columbia, MO 65211}

\email[P. G. Casazza]{pete@casazza.math.missouri.edu}
\email[M. C. Lammers]{lammers@sol.cms.uncwil.edu}

\subjclass{Primary: 46B15, 46B07}

\thanks{The first author was supported by NSF DMS 9706108 }

\noindent
\begin{abstract}
 We show that the classification problem for genus~$n$ Banach spaces can
 be reduced to  the unconditionally primary case and that the critical
 case there is $n=2$. It is further shown that a genus~$n$ Banach space 
is unconditionally primary if and only if it contains a complemented 
subspace 
of genus~$(n-1)$. We begin the process of classifying the genus~2
 spaces by showing they have a strong decomposition property.
\end{abstract}
  
\maketitle

\bigskip
\section{Introduction} \label{S:intro}

\TDS
It is well-known that a Banach space with a  basis has uncountably
many  non-equivalent normalized bases \ci{ps}.  However, there are
spaces with normalized unconditional bases that are unique up to
equivalence. G.~K\"othe and O.~Toeplitz \ci{kt} showed that $\ell_2$
has a unique unconditional basis and two papers by Lindenstrauss and
Pelczynski \ci{lp} and Lindenstrauss and Zippin \ci{lz} showed that the
complete list of spaces with a unique unconditional bases is $\ell_1,
\, \ell_2 $ and $c_0$.

One quickly notices that a unique normalized unconditional basis must
be symmetric.  This leads us to explore uniqueness up to a
permutation.  That is, two normalized unconditional bases are said to
be equivalent  up to a permutation if there exists a permutation of one
which is equivalent to the other.   Since the list of normalized
unconditional  bases that are actually unique  is now complete we will
use the phrase unique unconditional basis for unique up to a
permutation.  Edelstein and Wojtaszczyk \ci{ew} showed that direct
sums of $\ell_1, \, \ell_2 $ and $c_0$ have unique unconditional bases
and in 1985 Bourgain, Casazza, Lindenstrauss, and Tzafriri \ci{bclt}
showed that  2-convexified Tsirelson $T^2$ and  $(\sum _{n=1}^\infty
\oplus \ell_1)_{c_0}$,  $(\sum _{n=1}^\infty \oplus \ell_2)_{c_0}$,
$(\sum _{n=1}^\infty \oplus \ell_2)_{\ell_1}$, $(\sum _{n=1}^\infty
\oplus c_0)_{\ell_1}$, (along with  their complemented subspaces with
unconditional bases)   all have unique unconditional bases, while
somewhat surprisingly  $(\sum _{n=1}^\infty \oplus \ell_1)_{\ell_2}$
and  $(\sum _{n=1}^\infty \oplus c_0)_{\ell_2}$ do not.  More recently
Casazza and Kalton \ci{ck}, \ci{ck2} showed that other  Tsirelson type
spaces,  certain Nakano spaces, and some $c_0$ sums of $\ell_{p_n}$
with $p_n \to 1$ have unique unconditional bases.

In \ci{bclt} they define a new class of Banach spaces.  A Banach space
$X$ is said to be of {\bf genus}~${\bf n}$ if it and all it's complemented
subspaces with unconditional bases have  a unique  normalized
unconditional basis and there are exactly $n$ different complemented
subspaces with unconditional  bases, up to isomorphism.  For example,
the space $\ell_1 \oplus \ell_2$ is a Banach space of genus 3, the
three different complemented subspaces being $\ell_1 \oplus \ell_2$,
$\ell_1$ and $\ell_2$.  It was  shown in the memoir~\ci{bclt} that
2-convexified Tsirelson space, $T^2$,  is a Banach space of infinite
genus. In fact $T^2$ has uncountably many  non-isomorphic complemented
subspaces with  unconditional bases and every unconditional basis of a
complemented subspace is unique. It is unknown if there are any Banach
spaces of countable genus, i.e. of genus~$\omega$.  Our results apply
to this case as well so we include it in this paper.

At this point even the genus~2 spaces are unclassified, although there
is a conjecture that they are precisely the ones we already know(see
Appendix 5).  In section~2 we show that the problem of classifying
genus~n spaces reduces to classifying the unconditionally primary
genus~n spaces. We then characterize the unconditionally primary genus~$n$
 spaces  as those which contain a complemented subspace of genus~$(n-1)$.
 This basically shows that the backbone of this
classification problem is really the genus~2 case.  

  In section 3  we show that all Banach spaces of finite genus have
  the property that any subsequence of the original basis must contain
  a further subsequence equivalent to the unit vector basis of $c_0$,
  $\ell_1$ or $\ell_2$. This is particularly important to the genus~2
  case for it classifies such spaces into these three cases.

In section 4 we first show  that if the only spaces of genus 2
containing $c_0$ are those  conjectured by the memoir \ci{bclt} (
$(\sum _{n=1}^\infty \oplus \ell_1^n)_{c_0}$,  $(\sum _{n=1}^\infty
\oplus \ell_2^n)_{c_0}$), then the only space of genus 2  containing
$\ell_1$ are the duals of these spaces.  In other words it is enough
only to consider the ``$c_0$ case''. The remaining part of this
section deals with decomposing genus 2 spaces containing $c_0$. This
decomposition relies heavily on a result of Wojtaszczyk so we give some
details of this for clarity.

Finally we end with an appendix of a conjectured list of all  genus~$n$
spaces for$1\leq n \leq 6$.  We divide the genus~$n$ spaces into those
which are unconditionally primary and those which are not.

\section {Reducing  genus~$n$ spaces }

We start with a simple observation.  If $(x_n)$ and $(y_n)$ are sequences in
 Banach spaces $X$ and $Y$ respectively, we write $(x_n) \sim (y_n)$ to 
mean that $Tx_n = y_n$ defines an isomorphism from $span[x_n]$ to $span[y_n]$.
  Also we write $(x_n) \sim_ \pi (y_n)$ if there is a permutation $\pi$ of
 the natural numbers  so that $(x_n) \sim (y_{\pi(n)})$ .

\begin{prop} \lb{dualres}  If $X$ has an unconditional basis and $X^*$ has a
unique normalized unconditional basis, then $X$ has a unique
normalized unconditional basis.\
\end{prop}

\begin{proof}  Since $X^*$  is separable, every unconditional basis for
$X$ is shrinking.  So if $(x_n, \, x_n^*)$ and $(y_n, \, y_n^*)$ are
normalized unconditional bases for $X$ then  $(x_n^*)$ and $(y_n^*)$
are  bounded unconditional bases for $X^*$.  Hence $(x_n^*) \sim_\pi
(y_n^*)$ and so $(x_n) \sim_\pi (y_n)$.
\end{proof}

We had to assume that $X$ has an unconditional basis above since
$\ell_1$ has preduals without  unconditional bases.  We get
immediately:

\begin{cor} \lb{ducor}  If $X$ has a unconditional basis and $X^*$ is 
genus n,  then $X$ is of genus $\leq n$, for all $n\leq \omega$.
\end{cor}

Recall that a Banach space $X$ is said to be {\bf primary} if whenever
$X \cong Y \oplus  Z$, then either $X \cong Y$ or $X \cong  Z$.  We
say that $X$ is {\bf unconditionally primary} if whenever $X \cong Y
\oplus Z$ and $Y$, $Z$ have unconditional bases, then either $X \cong
Y$ or $X \cong Z$.  Now we wish to  give a characterization of
primary genus~$n$ Banach spaces. For this we need a recent result of
 Kalton \ci{K}.

\begin{thm}  If $X$ is a Banach space with an unconditional basis and $X$ 
has only countably many non-isomorphic complemented subspaces with 
unconditional bases  then $X \cong X^2$.
\end{thm}

The impact of the theorem is clear, for it gives  us that genus~$n$
spaces are isomorphic to their squares.  We are now ready for our
characterization of unconditionally primary spaces of finite genus.    

\begin{prop} Let $X$ be a Banach space of genus~n then $X$is unconditionally
  primary  if and only if $X$ contains a complemented subspace of 
\hbox{genus~$(n-1)$.} 

\end{prop}

\begin{proof} 
  $\Rightarrow$  Let $Z$ be a complemented subspace of $X$ with a
unique unconditional basis of maximal genus~$m$ and  $m<n$ .  If $m<
(n-1)$  there exists a complemented subspace $Y$ of $X$  with an
unique unconditional basis not equivalent to a complemented subspace
of $Z$ and not isomorphic to $X$.  Using the  theorem of Kalton,  we have
  $Y \cong Y \oplus Y$ and  $Y$ is complemented in $X$.  It  follows that
 $ X \cong X \oplus Y$.  Similarly  $ X \cong X \oplus Y \cong X \oplus Y
\oplus Z $. So $ Y \oplus Z$ is complemented in $X$ and has
genus~$k\leq m$ = genus($Y$).  But if $k=m$ then $Y\oplus Z \cong Y$
contradicting our assumption that $Z$ does not imbed complementably
into $Y$. Hence,  $k>m$ . But $Z$ was the  maximal  complemented space
with genus~$m$ and  $m<n$  so $ Y \oplus Z$  must be of genus~$n$ and
hence isomorphic to $X$.  So  $X$ is not primary contradicting our
assumption.  Therefore we must have  $m=n-1$.

$\Leftarrow$  By way of contradiction.  Suppose $W$ is a complemented
 subspace of  $X$ and $W$ is genus~$(n-1)$.  Now let $ X \cong Y
 \oplus Z$, where neither $Y$ nor $Z$ is isomorphic to $X$.  Then $Y$
 and $Z$  are genus $m_1$ and genus $m_2$ respectively with \hbox{
 $m_1, \, m_2 <n$} . 

 Since there are $n$-distinct unique unconditional bases for complemented
 subspaces of $X$,$(n-1)$ of them must be in $W$ and the remaining one is
the basis for $X$.  Hence if $(y_n)$ is a basis for $Y$ then because
$Y$ has genus $m_1$ with $m_1< n$ it can not be a basis for X hence it must
be equivalent to a subsequence of the basis for $W$. Similarly for $Z$.

This implies that $Y$  and $Z$ are both complemented subspaces of $W$
with unique unconditional bases, hence $ Y \oplus Z$ is a complemented
subspace of \hbox{ $W \oplus W \cong W$} by the preceding theorem. It
follows that  $X$ is a complemented subspace of $W$.   This is clearly
a contradiction since $X$ is genus~$n$ and $W$ is genus~$(n-1)$.

\end{proof}

Now we do the reduction of classifying genus~$n$ Banach spaces to the
unconditionally primary case.

\begin{thm}  Every Banach space $X$ of genus $n$ can be decomposed into $X
\cong X_1 \oplus X_2 \oplus \ldots \oplus X_m$,  where each $X_i$ is
unconditionally primary.
\end{thm}

\begin{proof}  If X is unconditionally  primary, we are done.  Otherwise
 by definition $X \cong  Y \oplus Z$ where neither $Y$ nor $Z$
 is isomorphic to $X$ 
and both $Y$ and $Z$ have unique unconditional bases.  Also
$Y$ and $Z$ are genus~$< n$.   Now iterate this
process until it stops.
\end{proof}

Theorem 2.5 tells us that to classify all Banach spaces of genus~$n$ we
only need to classify the unconditionally primary Banach spaces of
genus~$n$. From the results of \ci{bclt}, we know that the following
spaces are unconditionally primary (See 5. Appendix) :
\[(\sum_n \oplus E_n)_{c_0}, \, (\sum_n \oplus F_n)_{\ell_1}, \]
where $E_n$, $1 \leq n < \infty$ is 
\[ E_n=\ell_1^n \, {\mb or }\, \ell_1 ,\, {\mb or } \,\, \ell_2^n \,  
{\mb or }\, \ell_2, \]
and $F_n$, $1 \leq n < \infty$ is 
\[ E_n=\ell_\infty^n\, {\mb or }\, c_0 ,\, {\mb or }\, \ell_2^n \,
 {\mb or }\, \ell_2 .\]
It is natural then to conjecture that such iterations are the only way
to produce unconditionally primary spaces.  So we end this section
with the following conjecture.

\begin{con}  $X$ is unconditionally primary and genus n if and only if 
there is an unconditionally primary space $Y$ of genus $<$ n with
unconditional basis $(y_{i})$ and one of the following holds:
\begin{flushleft}
(1) \[X \cong(\sum \oplus Y)_{c_0} {\mb{ or }} X \cong(\sum \oplus
Y)_{\ell_1}\]

(2)\[X \cong(\sum \oplus Y_n)_{c_0} {\mb{ or }} X \cong(\sum \oplus
Y_n)_{\ell_1}\] where $Y_n =span[y_1,y_2, \dots y_n]$.

\end{flushleft}
\end{con}

This gives an indication of the  role that Banach spaces of genus~2
may play in the bigger picture of classifying all Banach spaces of
finite genus.  Conjecture 2.6 would look much more tractible if a
conjecture of \ci{bclt} were true.  That is, in \ci{bclt}, it is asked
if $X$ having a unique normalized unconditional  basis implies that
$c_0(X)$  also has a unique normalized unconditional  basis?
Recently, Casazza and Kalton \ci{ck2} have shown that this is false by
showing that $c_0$ sums of the original Tsirelson space fails to have
a unique normalized unconditional basis while in an earlier paper
\ci{ck} they showed that Tsirelson's space and its dual do have unique
normalized unconditional bases.

 
\section{Genus~$n$ spaces contain $c_0$, $\ell_1$, or $\ell_2$} 
 
If we consider only spaces of finite genus we get the following result.

\begin{thm} \lb{thm1} If $X$ is finite genus~$n$ then every normalized
unconditional basis for a complemented subspace of $X$ has a
subsequence which is permutatively equivalent to the unit vector basis
of $c_0$, $\ell_1$, or $
 \ell_2$.  \end{thm}

 To prove the Theorem we need  three propositions from \ci{bclt}.  The
first gives a condition on an unconditional basis which implies the
unconditional basis has  a permutation which is subsymmetric.  Recall
that an unconditional basis $(x_n)$ is {\bf subsymmetric} if it is
equivalent to all its subsequences.

\begin{prop}  [\ci{bclt}, Proposition 6.2] \lb{prp2} Let $X$ be a
Banach space with an unconditional basis $(x_n)$.  Suppose that every
subsequence of $(x_n)$ contains a further subsequence  which is
permutatively equivalent to $(x_n)$.  Then there exists a permutation
$\pi$ of the integers such that ($x_{\pi (n)}$) is a subsymmetric
basis.
\end{prop}

The next two propositions generalize results on Banach spaces with
symmetric  bases to  Banach spaces which have subsymmetric
bases. Actually Proposition \ref{prp4} generalizes a result on
homogeneous bases (bases  which are equivalent to all of  their
normalized block bases) but it is a well known result of  Zippin that
bases with this property must be equivalent to the unit vector basis
of $c_0$ or $\ell_p$ for some $1 \leq p < \infty$ and therefore are
symmetric.

\begin{prop}   [\ci{bclt}, Proposition 6.3] \lb{prp3}  Let $X$ be a
Banach space with a subsymmetric basis ($x_n$).  Let
$(\sigma_j)_{j=1}^{\infty}$ be mutually disjoint subsets of the
integers so that max $(\sigma_j)$ $ < $ min $(\sigma_{j+1})$, for all
$j$.  If we let $U$=$span[u_j = \sum_{i\in \sigma_j}{x_i}]$, then $X
\oplus U$ is isomorphic to $X$.  \end{prop}

\begin{prop}   [\ci{bclt}, Proposition 6.4] \lb{prp4} Let $X$ be a
Banach space with a normalized unconditional basis $(x_n)$.  Suppose
that for every normalized block basis with constant coefficients
$(u_j)$,  there exists a permutation $\pi$ of the integers so that
$(u_{\pi(j)})$ is equivalent to \xn.  Then \xn is equivalent to the
unit vectors in $c_0$ or $\ell_p$ for $1\leq p < \infty$.  \end{prop}

Combining Proposition \ref{prp3} and Proposition \ref{prp4} one can 
obtain the 
following immediate corollary \ci{bclt}.

\begin{cor} \lb{cor1} If $X$ has a subsymmetric basis and a unique
unconditional basis up to permutation, then $X$ is isomorphic to 
$c_0$, $\ell_1$, or $ \ell_2$.
\end{cor}

\begin{proof}   If   $(u_j)$ is a constant coefficient block basis of
the subsymmetric basis \xn then  by Proposition \ref{prp3}
$ span[u_j] \oplus X$ is isomorphic to $X$.  By the uniqueness of
the unconditional basis for $X$ the basis {$((u_j),(x_n))$ }is
permutatively equivalent to \xn.  Hence by Proposition \ref{prp4}, \xn is
equivalent to the unit vector basis of $c_0$ or $\ell_p$.  However for
$p\neq$ 1,2,  $\ell_p$ does not have a unique unconditional basis up to
a permutation.  This implies that \xn must be equivalent to the unit
vector basis of $c_0$ ,$\ell_1$ or $\ell_2$.  \end{proof}

 In particular, if $X$ has a subsymmetric basis and is genus n, then
 $X$ is isomorphic to $c_0$, $\ell_1$, or $ \ell_2$.
Now we are ready for the proof of the main theorem in this section.

\begin{flushleft}{\em Proof of Theorem \ref{thm1}}  By Corollary
\ref{cor1} it is enough to show that every normalized unconditional
 basis \xn  for a complemented subspace of $X$ has a subsymmetric
 subsequence.  Since $X$ is genus~$n$ there are only $n$ different
 normalized unconditional bases for complemented subspaces of $X$.
 Let $(x_m^i)_{m=1}^\infty, \,1\leq i\leq n$ be a representative of
 each of these $n$ different bases.  \end{flushleft}

\begin{claim}  $(x_n)$ has a subsequence $(x_n(1))$ with the property that
 every  subsequence of $(x_n(1))$ has a further subsequence
permutatively equivalent to $(x_n(1))$.  \end{claim}

 \begin{proof}
Either $(x_n)$ has the required property, and we are done, or $(x_n)$
has a subsequence $(x_n^1)$ which has no  further subsequence
equivalent to $(x_n)$.  Now, either $(x_n^1)$ satisfies the claim or
it has a subsequence $(x_n^2)$ which contains  no further subsequence
equivalent to $(x_n^1)$.  Continuing we find either  a  sequence
satisfying the claim or we can find sequences
$((x_k^i)_{k=1}^\infty)_{i=0}^{n-1}$, where $x_k^0 =x_k$, satisfying:
\begin{flushleft}

(1) $( x_k^i)_{k=1}^\infty$ is a subsequence of $(
x_k^{i-1})_{k=1}^\infty$, for  all $1\leq i\leq n-1$ 

(2)   $( x_k^i)_{k=1}^\infty$ has no subsequence equivalent to  $(
x_k^{i-1})_{k=1}^\infty$, for  all $1\leq i\leq n-1$ 
\end{flushleft}

Since $X$ is genus~$n$, it follows that $((
x_k^i)_{k=1}^\infty)_{i=0}^{n-1}$ must have exhausted the list, up to
permutatitve equivalence, of all  unconditional bases for a
complemented subspace of $X$.  But then by $(2)$, every subsequence of
$( x_k^{n-1})_{k=1}^\infty$ is permutatively equivalent to  $(
x_k^{n-1})_{k=1}^\infty$.

\end{proof}

By Claim 3.6 and Proposition 3.2 $(x_n)$ has a subsequence with a
permutation, call it $(y_n)$, so that $y_n$ is subsymmetric.  By
Corallary 3.5, $span[y_n]$ is isomorphic to $c_0$, $\ell_1$ or
$\ell_2$.  \begin{flushright}$\square$ \end{flushright}

The above argument works for Banach space of genus~$\omega$ also.  To
do this we need another result of Kalton \ci{K}.

\begin{prop}  If $X$ has an unconditional basis and at most countably many
subsequences of this basis span non-isomorphic Banach spaces, then $X$
is isomorphic to its hyperplanes.  \end{prop}

\begin{cor}  If a Banach space $X$ is of genus~$\omega$, then every
  normalized unconditional basis for $X$ has a subsequence  
 permutatively  equivalent to 
the unit vector basis of   $c_0$, $\ell_1$ or $\ell_2$. 
\end{cor}

\begin{proof} We will just note the changes required in the argument of
 Corollary~3.5.  Actually it is only Claim~3.6 that needs to be
 altered since the rest of the proof works perfectly well in this
 case.  Let $(x_n)$  be the unique normalized unconditional basis for
 $X$ and $(y_i^k)_{k,i=1}^\infty$ be a complete list of unconditional
 bases for  complemented subspaces of $X$.  Now either $(x_n)$ has the
 required property or $(x_n)$ has a subsequence $(x_n^1)$ which has no
 further subsequence equivalent to $(y_n^1)$.  Continue as in the
 proof of Claim~3.6, only now the process does not stop so we
 construct infinitely many subsequences with the properties:

(1)  $(x_k^i)_{k=1}^\infty$ is a subsequence of
$(x_k^{i-1})_{k=1}^\infty$.

(2)  $(x_k^i)_{k=1}^\infty$ has no subsequence equivalent to
$(x_k^{i-1})_{k=1}^\infty$ or   $(y_k^i)_{k=1}^\infty$.

Now choose the diagonal elements from these subsequences,
$(x_k^k)=(z_k)$ to get a subsequence of $(x_n)$.  There must be an $i$
so that  $(y_k^i)_{k=1}^\infty \sim (z_k)_{k=1}^\infty$.  But
$(z_k)_{k=i+1}^\infty$ is a subsequence of $(x_k^i)_{k=1}^\infty$ which
has no subsequence equivalent to  $(y_k^i)_{k=1}^\infty$.  Finally,
applying Proposition~3.7 we get the contradiction that
$(z_k)_{k=1}^\infty$ is equivalent to  $(z_k)_{k=i+1}^\infty$ (since
they span isomorphic spaces which by definition have unique
unconditional bases) while one of these is permutatively equivalent to
$(y_k^i)_{k=1}^\infty$ and the other is not. \end{proof}

One should note that this does not say that every subsequence of the
unconditional basis of a genus~n space  must contain the same $\ell_p$
unit vector basis.  Clearly $\ell_1 \oplus c_0$ has subsqeuences of
the unconditional basis equivalent to both that of $c_0$ and $\ell_1$.
However  for  a genus~2 space this is precisely the case and hence we
can  classify the genus~2 spaces into the categories of containing
$c_0$, $\ell_1$, or $\ell_2$.  We consider two of these cases in the
next section.


\section{Genus 2 spaces containing $c_0$}

First we notice by duality that if we can classify all spaces of
genus~2 containing $c_0$, then we also get the desired result for
those spaces containing $\ell_1$.

\begin{prop}  The only genus 2 spaces containing complemented $c_0$ are

\[(\sum_{n=1}^{\infty }\oplus \ell_2^n)_{c_0} 
\mb{ and }
(\sum_{n=1}^{\infty} \oplus \ell_1^n)_{c_0}\] 
 
if and only if  the only spaces of genus
2 containing complemented $\ell_1$
 are

 \[(\sum_{n=1}^{\infty} \oplus \ell_2^n)_{\ell_1}
 \mb{ and }`
(\sum_{n=1}^{\infty} \oplus \ell_{\infty}^n)_{\ell_1}\]
 
\end{prop}

\begin{proof}

$\Leftarrow$ Follows from Corollary\ref{ducor} 

$\Rightarrow$  We use the results of  James \ci{james} that   an
unconditional a basis for a space $X$ is boundedly complete iff the
space does not contain $c_0$, and the dual result relating shrinking
bases and $\ell_1$.  Let \xn be a normalized unconditional basis for a
genus 2 space $X$ containing $\ell_1$ complemented and let $(x_n^*)$
be the associated biorthogonal functions of \xn.  Then  $c_0$ can not
embed into $X$ or $X$ would be at least genus 3.  Hence \xn  is
boundedly complete, and so if $Y$ =$span[x_n^*]$, then $Y^* \cong X$.

\begin{claim}$Y$ is genus 2.  \end{claim}

\begin {proof} Let $(y_n)$ be a normalized unconditional basis 
 for a complemented subspace of $Y$ and let 
$(y_n^*)$ be the associated biorthogonal functions.  Because
the space spanned by $(y_n)$ can not contain complemented  $\ell_1$, 
again because $X$ did not  contain complemented  $c_0$, $(y_n)$ is a
shrinking basis for the $span[y_n)]$.
Since $Y^* \cong X$ and the basis is  shrinking it follows that 
$(y_n^*)$ is a normalized unconditional basis for a complemented subspace
of $X$.  Therefore, either $(y_n^*) \sim_{\pi} (x_n)$  for some permutation
$\pi$ or $(y_n^*) \sim (e_n)_{\ell_1}$. So $(y_n) \sim_{\pi} (x_n^*)$  
  or $(y_n) \sim (e_n)_{c_0}$.  Hence  $Y$ is genus 2.
\end{proof}

Now since \xn has subsequences equivalent to $ (e_n)_{\ell_1}$, $c_0$
embeds into $Y$.  Hence $Y$ is a genus 2 space containing $c_0$.  So $Y$
is isomorphic to 
\[(\sum_{n=1}^{\infty }\oplus \ell_2^n)_{c_0} \; {\rm or} \;
(\sum_{n=1}^{\infty} \oplus \ell_1^n)_{c_0}\]
and  $X$ is isomorphic to one of the duals of these two.  That is 

\[(\sum_{n=1}^{\infty} \oplus \ell_2^n)_{\ell_1} \;
 {\rm or} \; (\sum_{n=1}^{\infty} \oplus \ell_{\infty}^n)_{\ell_1} \, .\]

\end{proof}

One of the major difficulties in working with the genus~$n$ spaces is
that although every subsequence of the basis is equivalent to one of $n$
specified unconditional bases, there is no uniform constant of
equivalence.  For example, in \[(\sum \oplus \ell_1^n) _{c_0}\] the
natural basis of $\ell_1^n \oplus c_0$ becomes ``badly'' equivalent to
the unit vector basis of $c_0$ as $n$ increases.
  Our next goal is to produce a uniform
constant in genus~2 spaces for subsequences of the unconditional basis
which are equivalent to the whole basis.   That is, there is  a
constant $K$ so that any subsequence of the original basis of a genus
2 space  $X$ that spans a space  isomorphic to the original space, is
$K$-uniformly equivalent  to the original basis. Then we will use this
constant to show that $X$ has  a UFDD (unconditional finite
dimensional decomposition) of a very strong form.

In order to do this we first need a theorem of Wojtaszczyk \ci{Woj}
which also appears implicitly in \ci{mit}.  We should mention that
Wojtaszczyk's theorem can be applied to any Banach space with an
unconditional  basis not just those of genus 2.  In fact  the theorem
was originally used for a result on quasi-Banach spaces.  Although we
will not reproduce  the proof of this theorem we do need to present
some of the terminology and results from bipartite graph theory in
order to state the stronger version of Wojtaszczyk's theorem which we
need and which he actually  proved.   The following can be found in
Wojtaszczyk's paper \ci{Woj} and in more detail in \ci{b}.

A bipartite graph $G$ consists of two disjoint sets $N$ and $M$, and any set
 $E(G)$ of unordered pairs from $N \cup M$  with the
property that one element in the ordered pair is from $N$ and one is from
$M$.  We denote $N \cup M$  by $V(G)$.  We  call the elements of $V(G)$ the
vertices of the graph while $E(G)$ is called the edge set of $G$.  
  A subset $A\,  \subset V(G)$ is called one sided if 
$A\,  \subset N$ or $A\,  \subset M$.    Let $A $ be a one sided subset of 
$V(G)$ we say $A$ is matchable if there exists a 1-1 map $\psi :A \to
\mb{V(G)} $ such that $(a, \, \psi(a)) \in $ E(G) for all $a \in A$ and
we call $\psi $ a matching of $A$.
  
We now give a version of the classical Schroder-Bernstein
theorem of set theory, which has been observed by Banach \ci{ban}, 
in the language of bipartite graph theory.  

\begin{thm}  Let $M$, $N$ and $E(G)$ form  a bipartite graph $G$. 
  If both $M$ and $N$ are matchable  then there exists a matching
 of $N$, $\psi$ such that   $\psi(N)=M$.
\end{thm}

Now we are ready to   state Wojtaszczyk's theorem  and sketch how it
is proved.  This will elucidate the quantitative estimates needed.
We change the statement slightly, for the original theorem was stated
for quasi-Banach spaces and here we are only concerned with Banach
spaces.

\begin{thm}  If $(x_n)_{n \in N}$ and $(y_m)_{m \in M}$ are normalized
1-unconditional bases for Banach spaces $X$ and $Y$, and each is
equivalent to a permutation of a subsequence of the other, (that is
$(x_n)_{n \in N} \sim (y_{\sigma(n)})_{n \in N}$ for a 1-1 map
$\sigma:N\ \to  M$  and $(y_m)_{m \in M} \sim (x_{\gamma(m)})_{m \in
M}$ for a 1-1 map $\gamma:M \to  N$)  then $(x_n)_{n \in N}$ and
$(y_m)_{m \in M}$   are permutatively equivalent to each other.
\end{thm}

One should note that although $\sigma \mb{ and } \gamma$ are 1-1 maps
they need not be onto,  while the conclusion of the theorem implies
that there exists a 1-1 and onto map for the equivalence of $(x_n)$
and $(y_n)$.

 Wojtaszczyk uses bipartite graph theory and the classical
 Schroder-Bernstein theorem to obtain his result.  In particular he 
creates a bipartite graph $G$ with $V(G)= N \cup M$ and $E(G)=
\{ (n,\sigma (n))\}_{n \in N} \cup \{(m,\gamma (m))\}_{m \in M}$.  Since 
both $M$ and $N$ are matchable there exists a 1-1  map
$\Psi:N\, _{\overrightarrow{{\rm onto}}}\,  M$, and  a
 partition of $N$, $N=N_1 \cup N_2$, and hence a partition of $ M$, $M=
 \Psi (N_1)\cup \Psi (N_2)$ so that $(x_n)_{n \in N_1}$ is 
 equivalent to $(y_m)_{m \in \Psi( N_1)}$ and  $(x_n)_{n \in N_2}$ is
 equivalent to $(y_m)_{m \in \Psi( N_2)}$.  In particular

\bd
\Psi(n) = \left \{ \begin{array}{ll}
                   \sigma (n) & \mb{if $n \in N_1$} \\
                   \gamma ^{-1} (n) & \mb{if $n \in N_2$} \\
			   \end{array}
               \right. 
\ed
and therefore $(x_n)_{n \in N}$ and $(y_m)_{m \in
M}$   are permutatively equivalent to each other.  

If we consider the constants 
of equivalence $K_1$ and $K_2$ such that  
$(x_n)_{n \in N} \sim_{K_1} (y_{\sigma(n)})_{n \in N}$ and 
$(y_m)_{m \in M} \sim_{K_2} (x_{\gamma(m)})_{m \in M}$
by the 1-unconditionality of the basis one
can obtain:

\bea
 \| \sum _{n \in N}a_n  x_n \|  &\leq & 
\| \sum  _{n \in N_1} a_{n} x_{n} \| + 
\| \sum  _{n \in N_2} a_{n} x_{n} \| \no \\
&\leq &  K_1 \| \sum  _{n \in N_1} a_{n} y_{\Psi(n)} \| + 
 K_2 \| \sum  _{n \in N_2} a_{n} y_{\Psi(n)} \| \no \\
&\leq&  K_1 \| \sum  _{n \in N} a_{n} y_{\Psi(n)} \| + 
 K_2 \| \sum  _{n \in N} a_{n} y_{\Psi(n)} \| \no \\
&\leq&  ( K_1+  K_2) \| \sum  _{n \in N} a_{n} y_{\Psi(n)} \| \no .
\eea
The inequality in the other direction can be produced in a similar way.
The arguments above yield the following.

\begin {thm} Let  $(x_n)_{n \in N}$ and $(y_m)_{m \in M}$ be normalized
1-unconditional bases for Banach spaces $X$ and $Y$. If  $K_1$ and $K_2$
are constants  such that 
\begin{flushleft}
(1) $(x_n)_{n \in N} \sim_{K_1} (y_{\sigma(n)})_{n \in N}$ and

(2) $(y_m)_{m \in M} \sim_{K_2} (x_{\gamma(m)})_{m \in M}$ 
\end{flushleft}
then $(x_n)_{n \in N}$ and $(y_m)_{m \in M}$   are
$(K_1+K_2)$-permutatively equivalent to each other.
\end{thm}

An immediate corollary to this theorem is the form we will need for
the proof of the theorem below.

\begin{cor} \lb{wcor} If \xn is a normalized 1-unconditional basis
 and \xni is
 a subsequence of \xn which has a further subsequence which
 is $K$-equivalent
 to a  permutation of \xn, then \xni is $K+1$-equivalent to
 a permutation of 
\xn~. 
 \end{cor}

We are now ready to present a decomposition theorem for Banach spaces
of  genus~2.  To do this we first  produce a  uniform constant for
subbases of the original basis that span a space isomorphic  to the
original space.  This relies heavily on the theorem above  and the
fact that  any subsequence of an  unconditional basis for  a genus 2
Banach space $X$ containing  $c_0$ which is not equivalent to the unit
vectors in $c_0$, must be a basis for a space isomorphic $X$.  This is
clear since  any subsequence of the basis spans a complemented
subspace with an unconditional basis.   We start by  producing the
necessary uniform equivalence constant.

\begin{thm} \lb{unif} Let $X$ be a Banach space of genus 2  
 containing $c_0$ and let \xn be a normalized 1-unconditional basis
  for $X$.  Then there exists a natural number  K  such that for any
  subsequence \xni of \xn which spans a space isomorphic to  $X$, \xni
  is  K-permutatively equivalent to $(x_n)$.  \end{thm}

\begin{proof}   Assume no such $K$ exists and then we will proceed by
induction on $K$ to get a contradiction.  By Corollary \ref{wcor} there
exists a subsequence $(x_n^1)$ of \xn such that  
\begin{flushleft}
(i)  $span[(x_n^1)]$ is isomorphic to $X$ and

(ii)  no further subsequence of  $(x_n^1)$ is 2-permutatively
equivalent ( i.e. $K=2$) to \xn .
\end{flushleft}
If such a sequence did not exist and all subsequences that  spanned
the  space had a further subsequence that  was 2-permutatively
equivalent to \xn ,  then by the corollary  all subsequences of this
type would be 3-permutatively equivalent to \xn.
Let $k_0$=2 and choose $k_1 > k_0$  such that
\[ \| \sum _ {n=1}^{k_1} x_n^1 \| \geq 2 \]
Now since $(x_n^1)$ is permutatively equivalent to \xn and we have
 assumed
 no $K$ exists satisfying the theorem, there exists
 $(x_n^2)$ a subsequence of  $(x_n^1)$ such that

\begin{flushleft}
(i) $span[(x_n^2)]$ is isomorphic to $X$ and

(ii)  no further subsequence of  $(x_n^2)$ is $( k_0 +k_1)^2$-permutatively
 equivalent  to \xn .
\end{flushleft}

Without loss of generality we may assume that the  support of $(x_n^2)$ $>$
$k_0 +k_1$.   Now choose $k_2$ so that $k_2 > k_1$   and 

\[ \| \sum _ {n=1}^{k_2} x_n^2 \| \geq 2^2 .\]

Proceeding in this manner we can generate subsequences $(x_n^i)$ of \xn
such that for all $i \in N$ 

\begin{flushleft}
(1) $(x_n^{i+1})$ is a subsequence of $(x_n^i)$ with support 
of $(x_n^{i+1})$
$>$  $k_0+k_1 +k_2 + ... +k_i$ =$K_i$,

(2) $\| \sum _ {n=1}^{k_i} x_n^i \| \geq 2^i $,

(3) no subsequence of $(x_n^i)$ is  $(K_{i-1}) ^2$- permutatively equivalent
to \xn and

(4) $ span[(x_n^i)]=X$.

\end{flushleft}
First  consider the sequence
 $(z_n)=((x_n^j)_{n=1}^{k_j})_{j=1}^{\infty}$.  This is a normalized
 basis for a complemented subspace of $X$ and because of (2) above it
 can not span a space isomorphic to $c_0$.  Since $X$ is genus 2
 $(z_n)$ must be a basis for the whole space.  Also by genus 2, we
 know that $((x_n^j)_{n=1}^{k_j})_{j=1}^{\infty}$ is $D$-permutatively
 equivalent to \xn for some $D$ and some permutation $\pi$.   Now
 consider  the sequence  $(y_n^i)$ =$(x_n^i)_{n=1}^{K_{i-1}}  \cup
 ((x_n^j)_{n=1}^{k_j})_{j=i+1}^\infty$ for each $i \in \mathbb N $.
 Again this is a basis for the whole space  and is a  subsequence of
 $(x_n^i)$ by (1).

If we let $(a_n)$ be any set of scalars we have:

\bea 
\| \sum _ {n=1}^{\infty} a_n x_{\pi(n)} \| &\leq&  
 \| \sum _ {n=1}^{K_{i-1}} a_n x_{\pi(n)} \| 
	+ \| \sum _ {n=K_{i-1} +1 }^{\infty} a_n x_{\pi(n)} \| \no\\
&\leq&    K_{i-1} \| \sum _ {n=1}^{K_{i-1}} a_n y_n^i \| 
	+ D \| \sum _ {n=K_{i-1}}^{\infty} a_n z_n \| \no \\
&\leq&    K_{i-1} \| \sum _ {n=1}^{\infty} a_n y_n^i \| 
	+ D \| \sum _ {n=1}^{\infty} a_n y_n^i \| \no \\
&\leq&  (K_{i-1}+ D) \|\sum _ {n=1}^{\infty} a_n y_n^i \| 
\eea
 This follows from the fact that  $(z_n)$  and  $(y_n^i)$
 are the same sequence for $n > K_{i-1}$.
  Similarly we get,

\bea
\|\sum _ {n=1}^{\infty} a_n y_n^i\| &=&  
	\|\sum _ {n=1}^{K_{i-1}} a_n x_n^i +
	\sum _ {n=K_{i-1}+1}^{\infty} a_n z_n\| \no \\
&\leq&  \|\sum _ {n=1}^{K_{i-1}} a_n x_n^i \|+
	\| \sum _ {n=K_{i-1}+1}^{\infty} a_n z_n\|  \no \\
&\leq&  K_{i-1} \|\sum _ {n=1}^{\infty} a_n x_{\pi(n)}\| +
	 D \| \sum _ {n=1}^{\infty} a_n x_{\pi(n)}\|  \no \\
&\leq& (K_{i-1} + D) \| \sum _ {n=1}^{\infty} a_n x_{\pi(n)}\| 
\eea

By the two inequalities above we have that $(y_n^i)$ is $(K_{i-1} + D)$
 equivalent to a permutation of \xn. 
 But   the fact that  $(y_n^i)$ is a subsequence of 
$(x_n^i)$ and  from  (3) above we have  

\[ (K_{i-1} + D)  > (K_{i-1})^2 .\]

This is a contradiction since $K_i$ goes to $\infty$ and $D$ is fixed.
 
\end{proof}

Below  we will refer to this constant as  $K_0$ for a given Banach
space of genus~2 with unconditional basis \xn.
Now we show that spaces of genus~2 which contain $c_0$ have a strong
decomposition property.

A sequence $(E_n)_{n=1}^\infty$ of finite dimensional subspaces of a
Banach space $X$ is called an {\bf unconditional finite dimensional
decomposition}  for $X$, denoted UFDD, if for each $x \in X$ there is
a unique choice of $x_n \in E_n$ so that 
\[x=\sum_{n=1}^\infty x_n\] and the series is unconditionally convergent in
 $X$.  In this case we write \[X \cong \sum_{n=1}^\infty  \oplus E_n. \]  
If the $E_n$'s  are not necessarily finite dimensional we call this an 
{\bf unconditional Schauder decomposition}.   If there is a $K \geq 1$ so that 
$E_n$ is $K$-isomorphic to $X$, for all $n=1, 2, 3, ...$, we say
 that $X$ has an {\bf unconditional  decomposition into copies of itself}.

\begin{thm}  If $X$ is a Banach space of genus 2  containing  $c_0$, 
then there exists $X \cong \sum\oplus E_n$, an unconditional finite
dimensional decomposition satisfying:

\begin{flushleft}
(1) Each $E_n$ is dimension n and has a unconditional basis
$(x_i^n)_{i=1}^n$,

(2)  $((x_i^n)_{i=1}^n)_{n=1}^\infty$ is an unconditional basis for
$X$,

(3)  There is a constant $K \geq 1$ so that for all $n<m$ ,
\[(x_i^n)_{i=1}^n \sim_K (x_i^m)_{i=1}^n\]

(4) For all $n_1<n_2 \ldots $, \[ X \cong_{K} \sum\oplus E_{n_i}. \]
\end{flushleft}
\end{thm}

\begin{proof} Let \xn be a normalized unconditional basis for $X$.  By
Theorem~\ref{unif}, \xn is $K_0$-permutatively equivalent to
$(x_n)_{n=m}^\infty$, for $m=1,2, \ldots $.  Then letting $m=2$ above,
we can choose $n_1,n_2 >1$ so that $(x_i)_{i=1}^2$ is $K_0$-equivalent
to a permutation of $\{ x_{n_1}, x_{n_2}\}$.  Continuing we can find
$n_3, n_4, n_5 > max\{n_1,n_2 \}$ so that $(x_i)_{i=1}^3$ is
$K_0$-equivalent to a permutation of $\{ x_{n_3},x_{n_4}, x_{n_5}\}$.
After taking permutations we have found
$((x_i^n)_{i=1}^n)_{n=1}^\infty$,  a permutation of a subsequence of
\xn so that (1) holds with $E_n=span[(x_i^n)_{i=1}^n]$.  Now (2) of
the theorem is immediate and (3) follows from the choice of
$(x_i^n)_{i=1}^n$.  Finally, since $X$ is not isomorphic to $c_0$, 
\[\sup_n \| \sum_{i=1}^n x_i \| = +\infty .\]
Therefore, for every $n_1<n_2 \ldots $,
\[\sup_{n_i} \| \sum_{j=1}^{n_i} x_j \| = +\infty \]
and so $Y=\sum \oplus E_{n_i}$ is not isomorphic to $c_0$.  Hence $Y
\cong X$.  But $Y$ is the span of a subsequence of the basis for $X$.
So again by Theorem \ref{unif}, $Y \cong_{K_0} X$.

\end{proof}

We have immediately,

\begin{cor}  If $X$ is a Banach space of genus~2 containing $c_0$ then
 $X$ has an unconditional Schauder decompositon into copies of itself.  
i.e. $X \cong \sum \oplus  X$.
\end{cor}

The following result of Kalton~\ci{K} is the corresponding result for
  genus~$\omega$  Banach spaces.

\begin{prop}  If $X$ has an unconditional basis which has only countably many 
non-isomorphic subsequences, then $X$ has an unconditional Schauder
 decomposition into copies of itself.
\end{prop}

The next step to classifying genus~2 spaces would be to use the strong
decomposition result in Theorem~4.8 to show that genus~2 spaces
containing $c_0$ must be of the form $X=(\sum \oplus E_n)_{c_0}$.
After this there are enough tools available to complete the
classification. 
\section{Appendix} \label{S:append}
 
The following is a conjectured list of all genus~$n$ spaces for $1 \leq
n \leq 6$. 
 
\begin{flushleft}Genus 1: \end{flushleft}  The following are known \ci{bclt}
 to be the only genus~1 spaces:
\[ c_0,\, \ell_1, \, \ell_2 .\]
\begin{flushleft}Genus 2: \end{flushleft} The following are known  \ci{bclt}
to be genus~2 and have been conjectured \ci{bclt} to be the only spaces of 
genus~2. They are all unconditionally primary.
\[ (\sum_{n=1}^{\infty }\oplus \ell_2^n)_{c_0}, \,
(\sum_{n=1}^{\infty} \oplus \ell_1^n)_{c_0}, \, (\sum_{n=1}^{\infty}
\oplus \ell_2^n)_{\ell_1},  \, (\sum_{n=1}^{\infty} \oplus
\ell_{\infty}^n)_{\ell_1}  .\]
\begin{flushleft}Genus 3:\end{flushleft} The following are known  \ci{ew} to
 be genus~3: 
\[c_0 \oplus \ell_1, \, c_0 \oplus \ell_2, \, \ell_1 \oplus \ell_2 . \]
We conjecture the only unconditionally primary spaces of genus~3 are 
\[ (\sum _{k=1}^{\infty }\oplus (\sum _{n=1}^{\infty }
\oplus \ell_{\infty}^n)_{\ell_1}^k)_{c_0}, \,   (\sum_{k=1}^{\infty }
 \oplus (\sum_{n=1}^{\infty }
 \oplus\ell_1^n)_{\ell_\infty}^k)_{\ell_1} \]
\begin{flushleft}Genus 4:\end{flushleft}  We conjecture that the following
 are the only genus~4 spaces that are not unconditionally primary:
\bea  &(\sum _{n=1}^{\infty }\ell_2^n)_{\ell_1} \oplus (\sum
_{n=1}^{\infty }\ell_\infty^n)_{\ell_1},  \,  (\sum _{n=1}^{\infty
}\ell_2^n)_{\ell_1} \oplus  (\sum _{n=1}^{\infty
}\ell_\infty^n)_{\ell_1}, \no \eea and the following are the only
genus~4 spaces of that are unconditionally primary:  \bea &(\sum
_{k=1}^{\infty } \oplus(\sum _{j=1}^{\infty }\oplus  (\sum
_{n=1}^{\infty } \oplus \ell_1^n)_{\ell_\infty}^j)^k_{\ell_1})_{c_0}
\no  \\ &(\sum _{k=1}^{\infty } \oplus(\sum _{j=1}^{\infty }\oplus
(\sum _{n=1}^{\infty } \oplus
\ell_{\infty}^n)_{\ell_1}^j)^k_{\ell_\infty})_{\ell_1} \no   \eea
\begin{flushleft}Genus 5: \end{flushleft} We conjecture that the following 
 are   the only spaces of genus~5 that are not unconditionally primary:
\bea  & \ell_2 \oplus (\sum_{n=1}^{\infty} \oplus \ell_1^n)_{c_0}, \,
c_0 \oplus (\sum_{n=1}^{\infty} \oplus \ell_{\infty}^n)_{\ell_1} \no
\\ & \ell_1 \oplus (\sum_{n=1}^{\infty} \oplus \ell_1^n)_{c_0}, \,
\ell_2 \oplus (\sum_{n=1}^{\infty} \oplus \ell_{\infty}^n)_{\ell_1}
\no  \\ & \ell_2 \oplus (\sum_{n=1}^{\infty} \oplus \ell_2^n)_{c_0} ,
\, c_0 \oplus (\sum_{n=1}^{\infty} \oplus \ell_2^n)_{\ell_1} \no  \\ &
\ell_1 \oplus (\sum_{n=1}^{\infty} \oplus \ell_2^n)_{c_0}, \, \ell_2
\oplus (\sum_{n=1}^{\infty} \oplus \ell_2^n)_{\ell_1} ,  \no \eea and
the following are the only spaces of genus~5 that are unconditionally
primary: \bea & (\sum _{l=1}^{\infty }\oplus \sum _{k=1}^{\infty }
\oplus(\sum _{j=1}^{\infty }\oplus (\sum _{n=1}^{\infty } \oplus
\ell_1^n)_{\ell_\infty}^j)^k_{\ell_1})_{c_0}^l)_{\ell_1} \no  \\ &
(\sum _{l=1}^{\infty }\oplus \sum _{k=1}^{\infty } \oplus(\sum
_{j=1}^{\infty }\oplus (\sum _{n=1}^{\infty } \oplus
\ell_{\infty}^n)_{\ell_1}^j)^k_{\ell_\infty})_{\ell_1}^l) _{c_0} \no
\\  & (\sum _{k=1}^{\infty }\oplus (\sum _{n=1}^{\infty } \oplus
\ell_{2}^n)_{\ell_1}^k)_{c_0}, \,   (\sum_{k=1}^{\infty } \oplus
(\sum_{n=1}^{\infty } \oplus\ell_2^n)_{\ell_\infty}^k)_{\ell_1} . \no
\eea
\begin{flushleft}Genus 6:\end{flushleft}  We conjecture the following are the 
only spaces of genus~6 that are not unconditionally primary: \bea  &
 (\sum_{n=1}^{\infty }\oplus \ell_2^n)_{c_0} \oplus (\sum
 _{k=1}^{\infty } \oplus(\sum _{j=1}^{\infty }\oplus  (\sum
 _{n=1}^{\infty } \oplus \ell_1^n)_{\ell_\infty}^j)^k_{\ell_1})_{c_0}
 \no \\ & (\sum_{n=1}^{\infty} \oplus \ell_2^n)_{\ell_1} \oplus (\sum
 _{k=1}^{\infty } \oplus(\sum _{j=1}^{\infty }\oplus  (\sum
 _{n=1}^{\infty } \oplus
 \ell_{\infty}^n)_{\ell_1}^j)^k_{\ell_\infty})_{\ell_1}, \no \eea and
 that the following are the only spaces of genus~6 that are
 unconditionally primary: \bea & (\sum \oplus \ell_2)_{\ell_1}, \,
 (\sum \oplus c_0)_{\ell_1}, \, (\sum \oplus \ell_2)_{c_0}, \,  (\sum
 \oplus \ell_1)_{c_0} \no \\ & (\sum _{m=1}^{\infty }\oplus (\sum
 _{l=1}^{\infty }\oplus \sum _{k=1}^{\infty } \oplus(\sum
 _{j=1}^{\infty }\oplus (\sum _{n=1}^{\infty } \oplus
 \ell_1^n)_{\ell_\infty}^j)^k_{\ell_1})_{c_0}^l)_{\ell_1}^m)_{c_0} \no
 \, \\  & (\sum _{m=1}^{\infty }\oplus (\sum _{l=1}^{\infty }\oplus
 \sum _{k=1}^{\infty } \oplus(\sum _{j=1}^{\infty }\oplus (\sum
 _{n=1}^{\infty } \oplus
 \ell_{\infty}^n)_{\ell_1}^j)^k_{\ell_\infty})_{\ell_1}^l)
 _{c_0}^m)_{\ell_1}  \no . \eea


\end{document}